\documentclass{article}
\usepackage{amsmath,amssymb,latexsym,enumerate,amsthm}
\usepackage[dvips]{graphicx}
\newtheorem{Thm}[equation]{Theorem}
\newtheorem{Prop}[equation]{Proposition}
\newtheorem{Lem}[equation]{Lemma}

\theoremstyle{remark}

\newtheorem*{Rem*}{Remark}
\theoremstyle{definition}

\newtheorem{Ex}[equation]{Example}
\newtheorem*{Not*}{Notation}
\newtheorem{Def}[equation]{Definition}
\numberwithin{equation}{section}

\begin{document}

\title{Invariant group orderings and Galois conjugates}

\author{Peter A. Linnell\thanks{email:
linnell@math.vt.edu, url:
http://www.math.vt.edu/people/plinnell/, corresponding author}\\
Department of Mathematics\\
Virginia Tech\\
Blacksburg\\
VA 24061-0123\\
USA
\and
Akbar H. Rhemtulla\thanks{email:
akbar@math.ualberta.ca,
url: http://www.math.ualberta.ca/Rhemtulla\_A.html}\\
Department of Mathematics \\
University of Alberta \\
Edmonton \\
AL Canada T6G 2G1
\and
Dale P. O. Rolfsen\thanks{email:
rolfsen@math.ubc.ca,
url: http://www.math.ubc.ca/$\sim$rolfsen/}\\
Department of Mathematics \\
University of British Columbia \\
Vancouver \\
BC Canada V6T 1Z2}

\date{%
Tue Mar  4 13:07:29 EST 2008}
\maketitle

\begin{abstract}
This paper investigates conditions under which a given automorphism
of a residually torsion-free nilpotent group respects some ordering
of the group.  For free
groups and surface groups, this has relevance to ordering the
fundamental groups of three-dimensional manifolds which fibre over
the circle.



\end{abstract}

\section{Introduction}

We consider orderings of a residually torsion-free nilpotent
group $H$, that is a group $H$ which has a descending sequence
of normal subgroups $H_1 \supseteq H_2 \supseteq \cdots$
such that $H/H_i$ is torsion-free nilpotent for all $i$ and $\bigcap
H_i = 1$, and their invariance under
automorphisms, motivated by some questions in topology.  A group
$G$ is said to be \emph{ordered} if there is a strict total ordering
$<$ of its elements which is invariant under multiplication on
both sides, that is $f < g$ implies $hf < hg$ and $fh<gh$ for all
$f,g,h \in G$; we shall sometimes say that $G$ is \emph{bi-ordered}
in this situation, to emphasize the two-sidedness of the order.
The set $P$ of all $g \in G$ greater than the
identity is called the \emph{positive cone} of the ordering,
and satisfies:
\begin{enumerate}[\normalfont(1)]
\item \label{Iorder1}
$P$ is a sub-semigroup, that is, closed under multiplication, and $1
\not\in P$.

\item \label{Iorder2}
If $1 \ne g \in G$, then either $g \in P$ or $g^{-1}
\in P$, but not both.

\item \label{Iorder3}
If $g \in G$ and $p \in P$, then $g^{-1}pg \in P$; that is, $P$ is
normal in $G$.
\end{enumerate}
Conversely, if $P$ is a subset satisfying \eqref{Iorder1},
\eqref{Iorder2} and \eqref{Iorder3}, then it defines an ordering by
the formula $f < g \Leftrightarrow f^{-1}g \in P$.  If $P$ satisfies
only \eqref{Iorder1} and \eqref{Iorder2}, then it defines a
\emph{left-invariant} ordering and we say that $G$ is a
\emph{left-ordered} group.  Clearly an ordered group is left ordered,
and a left-ordered group is torsion free.

Many groups are orderable, including free groups and torsion-free
nilpotent groups, and more generally residually torsion-free
nilpotent groups.  Orderable groups have unique
roots: for $0 \ne n \in \mathbb{Z}$, we have $g^n = h^n$ if and only
if $g=h$ in $G$.  Another pleasant property of orderable groups is
that they obey the zero divisor conjecture: if $G$ is orderable and
$k$ is any integral domain, then the group ring $kG$ has no
nontrivial zero-divisors, and in fact embeds in a skew field.

A subset $X$ of a left-ordered group $G$ is \emph{convex} if
whenever $x,y,z \in G$ satisfy $x,z \in X$
and $x< y< z$, then $y \in X$.  The
collection of convex subgroups of a left-ordered group is linearly
ordered by inclusion.  An ordered group is \emph{Archimedean} if the
powers of every non-identity element are cofinal in the ordering.  By
theorems of H\"older and Conrad
\cite[Theorems 1.3.4 and 7.2.1]{MuraRhemtulla77}, every
Archimedean left-ordered group embeds (by a homomorphism preserving
the order) into the additive real numbers $\mathbb{R}$.

We say that an automorphism $\varphi \colon G \to G$ respects an
ordering $<$ if $g<h \Leftrightarrow \varphi(g) < \varphi(h)$.  In
this setting we also say that $<$ is a $\varphi$-invariant ordering.
This is equivalent to the equation $\varphi(P) = P$.

If $1 \to A \to B \to C \to 1$ is an exact sequence of groups and $A$
and $C$ are left-ordered, then $B$ can also be left-ordered by the
lexicographic order; specifically if we view $C$ as $B/A$, then for
$g,h \in B$, we define
$g>h$ if and only if $h^{-1}gA > A$, or $h^{-1}g \in A$ and $h^{-1}g
> 1$.  In this situation, the positive cone consists of the inverse
image of the positive cone of $C$ in $B$ together with the positive
cone of $A$.  This inheritance under extensions does not hold in
general for two-sided orderings; for example the Klein bottle
group, $\langle x,y \mid x^2=y^2\rangle$, is an extension of
$\mathbb{Z}$ by $\mathbb{Z}$ which is left orderable but
not orderable.  However the above recipe does provide a two-sided
ordering, provided the ordering of $A$ is respected by the
automorphisms of $A$ induced by conjugation by elements of $B$.
In particular, this holds if $A$ is central in $B$.

An application of invariant ordering is to the fundamental groups of
manifolds which fibre over the circle.  If $M^n$ fibres over $S^1$,
with fibre $F^{n-1}$, their groups fit into an exact sequence
\[
1 \longrightarrow \pi_1(F) \longrightarrow \pi_1(M)
\longrightarrow \pi_1(S^1) = \mathbb{Z} \longrightarrow 1.
\]
One may consider $M$ as a mapping torus of a
monodromy map $f \colon F \to F$, that is
\[
M \cong F \times [0,1]/(x,1) \sim (f(x),0),\quad x \in F.
\]
The fundamental group $\pi_1(M)$ is an HNN extension of $\pi_1(F)$
defined by the automorphism $f_* \colon \pi_1(F) \to \pi_1(F)$, and
so it is easy to see the following.
\begin{Prop} \label{Ppiorderable}
If $M$ is a manifold which fibres over the circle with fibre $F$,
then $\pi_1(M)$ is left orderable if and only if $\pi_1(F)$ is
left orderable.  Moreover $\pi_1(M)$ is orderable if and only if
$\pi_1(F)$ has a (two-sided) ordering invariant under the
$\pi_1$-monodromy $f_* \colon \pi_1(F) \to \pi_1(F)$.
\end{Prop}

We will now consider lower dimensions.  By a surface, we understand a
metric space each of whose points has a neighborhood homeomorphic with
the Euclidean plane or half-plane.  It is known that almost all
surface groups are orderable; more precisely \cite[Theorem
3]{RolfsenWiest01} yields
\begin{Prop} \label{Psurface}
If $F$ is any connected surface (compact or not), then $\pi_1(F)$
is orderable unless $F$ is a projective plane or Klein bottle.  The
Klein bottle group is left orderable but not orderable.
\end{Prop}
In particular, all orientable surface groups are orderable.  The
nonorientable closed surfaces have orderable groups if and only if
the surface has negative Euler characteristic, or equivalently, the
surface is the connected sum of at least three projective planes.
The projective plane's fundamental group, being finite, is certainly
not left orderable.

The following theorem is \cite[Theorem 1.1]{PerronRolfsen03} in the
case that $\pi_1(F)$ is free, as in fibred knot complements $M = S^3
\setminus N(K)$.  The more general case in which $F$ may be a
\emph{closed} orientable surface was proved in
\cite[Corollary 2.3]{PerronRolfsen06}.
\begin{Thm} \label{TPR}
Suppose that $M^3$ is an orientable 3-manifold which fibres over
$S^1$, with compact orientable fibre $F^2$ and monodromy $f \colon F
\to F$.  Then $\pi_1(M)$ is orderable if all the eigenvalues of the
homology monodromy $f_* \colon H_1(F) \to H_1(F)$ are real and
positive.  In particular, a fibred knot in $S^3$ or a homology
sphere, has orderable knot group if all the roots of its Alexander
polynomial are real and positive.
\end{Thm}

One of the main points of the present paper is to investigate the
extent to which the condition on the eigenvalues is necessary in
Theorem \ref{TPR} above.  Our main result for producing examples is
Proposition \ref{Pgeneral}.

We would like to thank Thomas Schick for pointing out a couple of
inaccuracies in a preliminary version of this paper, and the referee
for some useful comments.

\section{Abelian orderings and Galois conjugates}

Let $\theta$ be an endomorphism of the finite rank torsion-free
abelian group
$A$.  The eigenvalues of $\theta$ will mean the complex eigenvalues
of the $\mathbb{C}$-linear transformation $\theta \otimes 1$ induced
by $\theta$ on the finite dimensional $\mathbb{C}$-vector space
$A \otimes_{\mathbb{Z}}
\mathbb{C}$.  If two algebraic complex numbers have the same minimal
polynomial over $\mathbb{Q}$, they are said to be Galois conjugates.

\begin{Lem} \label{Leigenvalues}
If the finite rank abelian group $A \ne 1$ is a
subgroup of the additive
group of real numbers $\mathbb{R}$ and the endomorphism $\theta
\colon A \to A$ is multiplication by the real number $\alpha$, then
the eigenvalues of $\theta$ are the Galois conjugates of $\alpha$.
\end{Lem}
\begin{proof}
Let $f$ denote the minimal polynomial of $\alpha$ over $\mathbb{Q}$.
Then $f(\alpha)=0$ and $f$ is irreducible in $\mathbb{Q}[X]$,
more or less by definition.
Observe that if $g \in \mathbb{Z}[X]$, then $g(\theta)$ is given by
multiplication by $g(\alpha)$.  It follows that $f(\theta) = 0$.
Therefore the eigenvalues of $\theta$ satisfy $f$.  Since $f$ is
irreducible, it follows that the eigenvalues of $\theta$ are
precisely the roots of $f$, in other words the eigenvalues of
$\theta$ are the Galois conjugates of $\alpha$.
\end{proof}

\begin{Prop} \label{Pmain}
Let $A$ be a torsion-free abelian group of finite rank and let
$\theta$ be an automorphism of $A$.
Then $\theta$ preserves an order if and only if for
each eigenvalue of $\theta$, at least one of its Galois conjugates
is a positive real number.
\end{Prop}
\begin{proof}
We may assume that $A\ne 1$.
First suppose $\theta$ preserves an order on $A$.
If $1 = A_0 < A_1 < \dots < A_n = A$ is a series of convex subgroups
in $A$, then each $A_{i+1}/A_i$ is a torsion-free abelian group with
rank at least 1, consequently $n$ is at most the rank of $A$.
Therefore we have a finite
series of convex subgroups $1 = A_0 < A_1 < \dots < A_n = A$ of $A$,
where $A_{i-1}\ne A_i$ and there are no convex subgroups strictly
between $A_{i-1}$ and $A_i$ for all $i$.  Since $\theta$ maps convex
subgroups to convex subgroups, it follows that $\theta A_i = A_i$ for
all $i$.  Thus $\theta$ induces an
order preserving automorphism on the ordered group $A_i/A_{i-1}$.  By
H\"older's theorem \cite[Theorem 1.3.4]{MuraRhemtulla77}
we may consider $A_i/A_{i-1}$ as a subgroup of the
additive group of $\mathbb{R}$ (with its natural order), and then
a theorem of Hion \cite[Theorem 1.5.1]{MuraRhemtulla77} tells us that
this automorphism induced by $\theta$ is multiplication by a
positive real number.  It follows
from Lemma~\ref{Leigenvalues} that for each
eigenvalue of $\theta$, at least one of its Galois conjugates is a
positive real number.

Conversely suppose all the eigenvalues of $\theta$
have a Galois conjugate which is a positive real number.
If $1 = A_0 < A_1 < \dots < A_n = A$ is a sequence of subgroups of
$A$ with each $A_{i+1}/A_i$ torsion-free, then $n$ is at most the
rank of $A$.  Therefore we may choose a finite series of
$\theta$-invariant subgroups $1 = A_0 < A_1 < \dots < A_n = A$ where
$A_i/A_{i-1}$ is torsion free and if $B$ is a $\theta$-invariant
subgroup such that $A_{i-1} < B \le A_i$, then $A_i/B$ is a
torsion group.

We now fix $i$ and view $A_i/A_{i-1}$ as a $\mathbb{Z}[X]$-module,
where the action of $X$ is induced by $\theta$.  Then $A_i/A_{i-1}
\otimes_\mathbb{Z} \mathbb{Q}$ is an irreducible
$\mathbb{Q}[X]$-module.
By the structure theorem for modules over a principal ideal domain
applied to the ring $\mathbb{Q}[X]$ and the $\mathbb {Q}[X]$-module
$A_i/A_{i-1} \otimes_{\mathbb{Z}} \mathbb{Q}$, we may consider
$A_i/A_{i-1}$ as a subgroup of $\mathbb {R}$, and $\theta$ as
multiplication by a positive real number, because each eigenvalue of
$\theta$ has at least one of its Galois conjugates a positive real
number.  Thus each $A_i/A_{i-1}$ has a $\theta$-invariant order.
We can now use the lexicographic ordering to obtain a
$\theta$-invariant order of $A$.
\end{proof}

\section{Orderings of Residually Torsion-free Nilpotent Groups}

Suppose $F$ is a finitely generated nonabelian free group and
$\theta \colon F \to F$ is an
automorphism.  It was shown in \cite[Theorem 2.6]{PerronRolfsen03}
that if (as in
Theorem \ref{TPR}) all the eigenvalues of the corresponding map on
the abelianization $H_1(F,\mathbb{Z})$ are real and positive, then
one can construct a $\theta$-invariant order on $F$.  In this
section we show that the condition that all the eigenvalues be
positive is not necessary.  We will also consider generalizations to
residually torsion-free nilpotent groups.

We need to consider the rational lower central series $G_n^r$ (for
$n$ a positive integer) of the group $G$; here the superscript $r$
indicates ``rational".  Recall
that the lower central series of a group
$G_n$ for a group $G$ is defined inductively
by $G_1 = G$ and $G_{n+1} = [G,G_n]$ for $n \ge 1$.
Then $G_n^r = \{g \in G \mid g^m \in G_n$ for some positive integer
$m\}$.  Since $G/G_n$ is a nilpotent group, $G_n^r$ is
a characteristic subgroup of $G$ for all positive integers $n$.
Furthermore, each $G_n^r/G_{n+1}^r$ is a torsion-free group which
lies in the center of $G/G_{n+1}^r$.  It is easy to see that
$\bigcap_{n=1}^{\infty} G_n^r = 1$ if and only if
$G$ is residually torsion-free nilpotent.  If $\theta$ is
an automorphism of $G$, then $\theta$ induces automorphisms
on $G_n/G_{n+1}$ and $G_n^r/G_{n+1}^r$ for all $n \ge 1$.
Thus if $A$ is the
infinite cyclic group, we can make each $G_n/G_{n+1}$ into a
$\mathbb{Z}A$-module by making the generator of $A$ act by $\theta$.
Then \cite[5.2.5]{Robinson96} tells us that we have a
$\mathbb{Z}A$-epimorphism from $(G/G_2)^{\otimes n}$, the $n$-fold
tensor product with diagonal $A$-action, onto $G_n/G_{n+1}$.
We may also view each $G_n^r/G_{n+1}^r$ as a $\mathbb{Z}A$-module.
It follows that we have a $\mathbb{Q}A$-epimorphism from
$(G/G_2^r \otimes_{\mathbb{Z}} \mathbb{Q})^{\otimes n}$ onto
$(G_n^r/G_{n+1}^r) \otimes_{\mathbb{Z}} \mathbb{Q}$.
Suppose now that $G/G_2^r$ has finite rank and
let $\alpha_1, \dots,\alpha_p$ denote the
eigenvalues of the automorphism induced by $\theta$ on $G/G_2^r$.
Then we see that $G_n^r/G_{n+1}^r$ is a
finite rank torsion-free abelian group for all $n$, and that the
eigenvalues of the automorphism induced by $\theta$
on each $G_n^r/G_{n+1}^r$
are just products of $\alpha_1,\dots,\alpha_p$.

To apply Proposition \ref{Pmain},
we want a condition that will ensure that at least one of the Galois
conjugates of such a product is a positive real number.
For convenience, we make the following definition.
\begin{Def} \label{Dspecial}
Let $f \in \mathbb{Q}[X]$ be a monic polynomial
and let $f = f_1\dots f_n$ be its factorization
into monic irreducible polynomials (so $n$ is a nonnegative
integer and each $f_i$ is irreducible).  Then we say that $f$ is a
\emph{special} polynomial if for each $i$, at least one of the
following two conditions is satisfied.
\begin{enumerate}[\normalfont(i)]
\item\label{Dspeciali}
$f_i$ has odd prime power degree, negative constant term, and all its
roots are real.

\item\label{Dspecialii}
All the roots of $f_i$ are real and positive.
\end{enumerate}
\end{Def}
In an earlier version of this paper, we defined $f$ to be special if
it always satisfied condition \ref{Dspecial}\eqref{Dspeciali} (a
stronger condition).  With this definition, Wangshan Lu in his
M.Sc.\ thesis \cite{Lu07}
classified all the special Alexander polynomials of
fibred knots with degree less than 10.

We can now state

\begin{Lem} \label{Lpositive}
Let $f \in \mathbb{Q}[X]$ be a special polynomial.
If $\alpha$ is a product of the roots of $f$,
then at least one of the Galois conjugates of
$\alpha$ is a positive real number.
\end{Lem}
\begin{proof}
Let $f = f_1\dots f_n$ be the factorization of $f$ into monic
irreducible polynomials.
We prove the result by induction on $n$, the case $n=0$ being clear
because $\alpha = 1$ (empty product) in this case.
Let $\alpha_1,\dots, \alpha_p$
be the roots of $f_1$, let $\beta_1, \dots, \beta_q$ be the roots
of $f_2 \dots f_n$, and let $K = \mathbb{Q}(\alpha_1,
\dots,\alpha_p,\beta_1, \dots,\beta_q)$.
Then
\[
\alpha = \alpha_1^{r_1} \dots \alpha_p^{r_p} \beta_1^{s_1} \dots
\beta_q^{s_q}
\]
where the $r_i,s_i$ are nonnegative integers.  Clearly all the
Galois conjugates of $\alpha$ are real numbers.
Set $\beta = \beta_1^{s_1} \dots \beta_q^{s_q}$.  Suppose all the
Galois conjugates of $\alpha$ are negative.  First consider the case
when all the roots of $f_1$ are positive.  By induction on $n$, some
Galois conjugate of $\beta$ is positive and it follows that some
Galois conjugate of $\alpha$ is also positive.

Therefore we may assume that $f_1$ satisfies condition
\ref{Dspecial}\eqref{Dspeciali}, and here we may assume that $p$ is
a power of the odd prime $p'$.  Let $c$ denote the
negative of the constant term of $f_1$, a positive real number.
By considering
a Sylow $p'$-subgroup of the Galois group of $K$ over $\mathbb{Q}$,
there is a $p'$-subgroup $P$ of field automorphisms of $K$ which acts
transitively on $\{\alpha_1, \dots, \alpha_p\}$.  Set $\bar{p} =
|P|/p$, so $\bar{p}$ is the order of the stabilizer in $P$ of
$\alpha_1$.  Then we have
\[
\prod_{\theta \in P} \theta(\alpha) = (\alpha_1 \dots
\alpha_p)^{\bar{p}(r_1 + \dots + r_p)} \prod_{\theta \in P}
\theta(\beta) = c^{\bar{p}(r_1 + \dots + r_p)}
\prod_{\theta \in P} \theta(\beta).
\]
Write $\gamma = \prod_{\theta \in P} \theta(\beta)$, so
$\prod_{\theta \in P} \theta(\alpha) =
c^{\bar{p}(r_1 + \dots + r_p)} \gamma$.
Since $\gamma$ is a product of the roots of $f_2 \dots f_n$,
by induction on $n$ there is a field automorphism $\phi$ of $K$ such
that $\phi(\gamma)$ is a positive real number.  Then
\begin{equation} \label{Econtradiction}
\prod_{\theta \in P} \phi\theta(\alpha) =
c^{\bar{p}(r_1 + \dots + r_p)} \phi(\gamma).
\end{equation}
The left hand side of \eqref{Econtradiction} is a product of $|P|$
negative numbers and hence is negative, whereas the right hand side
is positive.  This contradiction completes the proof.
\end{proof}

\begin{Prop} \label{Pgeneral}
Let $f \in \mathbb{Q}[X]$ be a special polynomial (see Definition
\ref{Dspecial}), let $G$ be a residually torsion-free nilpotent group,
let $\theta$ be an automorphism of $G$, and let $\phi \colon G/G_2^r
\to G/G_2^r$ be the automorphism induced by $\theta$.  Assume that
$G/G_2^r$ has finite rank and that the eigenvalues of $\phi$ are
roots of $f$.  Then $G$ has a bi-ordering invariant under $\theta$.
\end{Prop}
\begin{proof}
Recall that if $Z$ is a central subgroup of the group
$G$, and $Z$ and $G/Z$ are orderable, then one can use the
lexicographic ordering to bi-order $G$.
To show that $G$ has a bi-order
which is invariant under $\theta$, it will be sufficient to show that
each $G/G_n^r$ has a bi-order which is invariant under $\theta$, by
\cite[Theorem 1.3.2(a)]{MuraRhemtulla77}.
Thus it will be sufficient to show that each
$G_n^r/G_{n+1}^r$ has an order which is invariant under $\theta$.
Lemma \ref{Lpositive} tells us that each
eigenvalue of the induced action of $\theta$ on $G_n^r/G_{n+1}^r$ has
at least one of its Galois conjugates a positive real number.
We now see from Proposition \ref{Pmain} that $G_n^r/G_{n+1}^r$ has an
ordering invariant under $\theta$ and the result follows.
\end{proof}

\begin{Rem*}
To apply Proposition \ref{Pgeneral}, we need examples of residually
torsion-free nilpotent groups $G$ such that $G/G_2^r$ has finite
rank.  The latter condition is easily satisfied; it will certainly be
the case if $G$ is finitely generated.  The former is certainly
satisfied if $G$ is torsion-free nilpotent or free.

Recall that the group $G$ is \emph{fully residually free} means that
given $g_1,\dots g_n \in G$, then there exists a homomorphism $\theta
\colon G \to F$ where $F$ is a free group such that $\theta(g_i) \ne
1$ for all $i$.  Certainly a fully residually free group is
residually torsion-free nilpotent.
Benjamin Baumslag \cite[top of p.~414]{Baumslag67} proved that if
$4 \le n \in \mathbb{Z}$ and $0\ne w_1,\dots,w_n \in \mathbb{Z}$,
then $\langle a_1,\dots,a_n \mid a_1^{w_1}\dots a_n^{w_n} =
1\rangle$ is fully residually free.  Thus in
particular all surface groups, orientable or not, are residually
torsion-free nilpotent, except the projective plane, Klein bottle,
and the group $\langle a,b,c \mid a^2b^2c^2 = 1\rangle$, the
fundamental group of the connected sum of exactly three
projective planes (the torus group is fully residually free, though
this particular case does not follow from Baumslag's paper).
The group $\langle a,b,c \mid a^2b^2c^2 = 1\rangle$ is not fully
residually free \cite{Lyndon59}; we do not know whether or not it is
residually torsion-free nilpotent.

Another situation where Proposition \ref{Pgeneral} can be applied is
as follows.
In \cite[p.~17]{Labute70}, Labute defines an element $x$ of a free
group $F$ to be \emph{primitive} if $x \ne 1$, and if $x \in
F_n\setminus F_{n+1}$ (where $F_n$ denotes the lower central series
of $F$), then $x$ is not a $d$th power modulo $F_{n+1}$ for any
integer $d \ge 2$.  Then \cite[Theorem on p.~17]{Labute70}
shows that a
one-relator group, where the relator is a primitive element of the
ambient free group, has a lower central series with torsion-free
factors.  Thus in particular if the one-relator group is also
residually nilpotent, then the group is residually torsion-free
nilpotent.
\end{Rem*}

\begin{Ex}
Let $p$ be a power of an odd prime,
let $F$ denote the free group of rank $p$
with free generators $x_1, \dots, x_p$,
and let $f(X) = X^p + f_{p-1}X^{p-1} + \dots
+ f_1X -1 \in \mathbb {Z}[X]$ be a polynomial which is irreducible in
$\mathbb {Q}[X]$ and has all roots real; note that $f$ is a special
polynomial (Definition \ref{Dspecial}).  For example $p=3$ and $f(X)
= x^3 - 3x -1$, which has two negative roots.
Then we can define an automorphism $\theta$ of $F$
by the formula
\[
\theta x_1 = x_2,\ \theta x_2 = x_3,\ \dots,\ \theta x_{p-1} = x_p,\
\theta x_p = x_1 x_2^{-f_1}\dots x_p^{-f_{p-1}}.
\]
Then the eigenvalues of the automorphism of $F$ induced by $\theta$
are roots of $f$.  Since $F/F_2^r$ has finite rank, namely $p$, and
$F$ is residually torsion-free nilpotent, it follows from Proposition
\ref{Pmain} that $F$ has an ordering which is invariant under
$\theta$.
\end{Ex}

\begin{Ex} \label{Erelation}
Let $n$ be a positive integer, let $F
= \langle x,y\rangle$ be the free group of rank 2 with generators
$x,y$, and let $G = \langle g \mid g^n = 1\rangle$ be the cyclic
group of order $n$.  By mapping $F$ onto $G$ by the map $x \mapsto
g$, $y \mapsto 1$, we may write $G = F/R$ where $R$ is the kernel of
this map.  Then $R$ is a finitely generated free group and
conjugation by $F$ on $R$ gives $R/R'$ (where $R'$ indicates the
commutator subgroup of $R$) the structure of a
$\mathbb{Z}G$-module.  Furthermore,
\cite[Proposition 5.10]{Gruenberg76} shows that $R/R' \cong
\mathbb{Z} \oplus \mathbb{Z}G$.  Let $\theta$ indicate the
automorphism $r \mapsto xrx^{-1}$ of $R$ and let $\phi$ be the
automorphism of $R/R'$ induced by $\theta$.  Since $\phi$ has order
$n$, the eigenvalues of $\phi$ are precisely $n$th roots of 1.
Also by restricting a bi-order on $F$ to $R$, we see that $\theta$
preserves a bi-order on $R$.
\end{Ex}

This contrasts with the example of an automorphism of a free group
which is periodic of period $n$, where $2 \le n \in \mathbb{Z}$; that
is $\theta^n$ is the identity, but $\theta$ is not the identity.  Its
eigenvalues are also $n$th roots of 1.  However no ordering of the
free group can be $\theta$-invariant.  For if (say) $w \in F$ and $w
< \theta(w)$, then $\theta(w) < \theta^2(w)$, etc.\ and we obtain the
contradiction $w < \theta^n(w) = w$.

This shows in particular, that one cannot determine whether an
automorphism of a finitely generated free group preserves a bi-order
by looking at its action on the abelianization of the free group.

Finally we give an explicit example of a fibred knot whose group
is ordered, but is not covered by the criterion of \cite[Theorem
1.1]{PerronRolfsen03}, that is a fibred knot whose Alexander
polynomial doesn't have all roots real and positive.
Consider the polynomial $f(x) :=x^6+3x^5-x^4-7x^3-x^2+3x+1$.
This factors
as $(x^3+x^2-2x-1)(x^3+2x^2-x-1)$.  It is easily checked that both
factors are irreducible and have all roots real, hence $f(x)$ is a
special polynomial (Definition \ref{Dspecial}).  Furthermore $f$ has
negative roots.  If we multiply
$f(x)$ by $-x^{-3}$, we obtain $-x^3-3x^2+x+7+x^{-1}-3x^{-2}-x^{-3}$.
This is symmetric in $x$ and $x^{-1}$, the sum of its
coefficients is 1, and its leading coefficient is $-1$,
hence it is the Alexander polynomial of a fibred
knot.  The corresponding Conway polynomial is $\nabla(z) :=
1 -20z^2 - 9z^4 -z^6$, because $\nabla(x^{-1/2} +x^{1/2}) =
-x^3-3x^2+x+7+x^{-1}-3x^{-2}-x^{-3}$.
An explicit example of a fibred knot with Conway polynomial
$\nabla$ is given below.
\smallskip

\includegraphics[scale=.6]{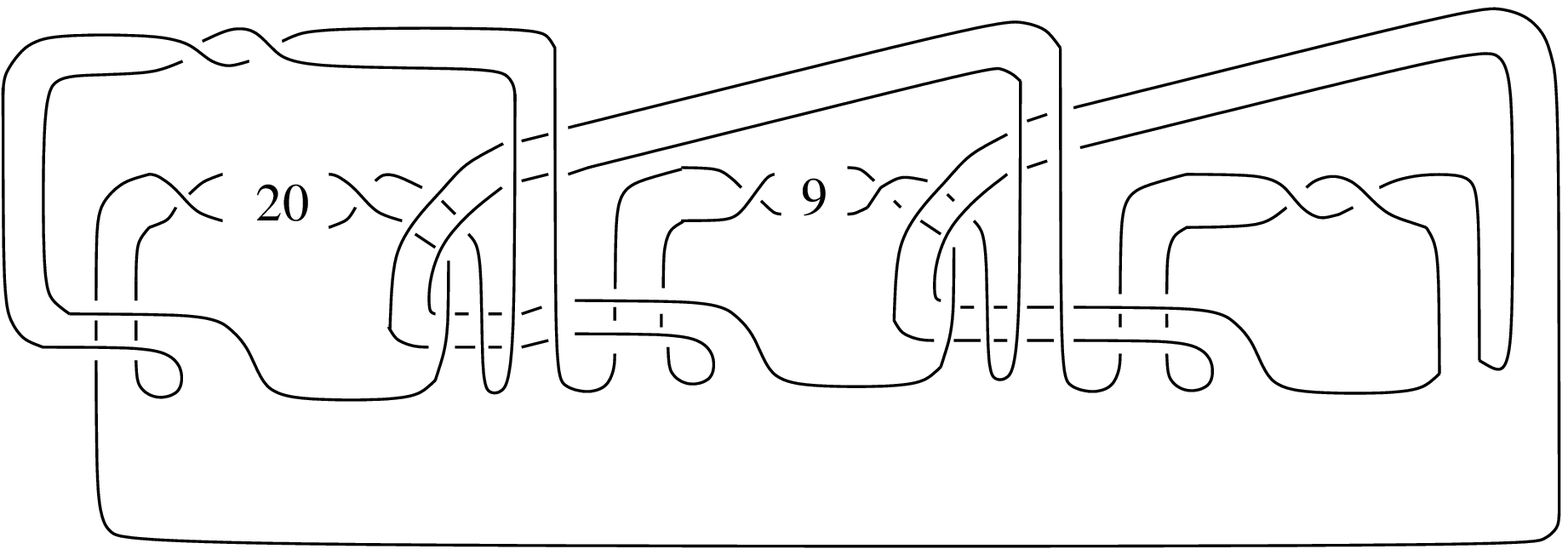}

This follows from \cite[Theorem 6]{Morton83}; in that theorem the
Conway polynomial is called the Kauffman polynomial.
The numbers 20 and 9 indicate the number of full twists in the
direction indicated (a total of 40 and 18 crossovers, including
the ones drawn).

\bibliographystyle{plain}

\end{document}